\documentclass[12pt]{elsarticle}

\usepackage{latexsym,amsfonts,amsmath,theorem,amssymb}
\usepackage{graphicx}

\def\be{\begin{equation}}
\def\bea{\begin{eqnarray*}}
\def\ee{\end{equation}}
\def\eea{\end{eqnarray*}}
\def\ba{\begin{array}}
\def\ea{\end{array}}
\def\bi{\begin{itemize}}
\def\ei{\end{itemize}}
\newtheorem{theo}{Theorem}
\newtheorem{lem}{Lemma}

\def\tP{\tilde{P}}
\def\tG{\tilde{G}}
\def\tT{\tilde{T}}
\def\tf{\tilde{f}}
\def\tg{\tilde{g}}

\journal{XXX}

\begin{document}
\begin{frontmatter}

\title{Haar system as Schauder basis in Besov spaces: The limiting cases
for $0<p\le 1$}
\author[bonn]{Peter Oswald}
\ead{agp.oswald@gmail.com}

\address[bonn]{Institute for Numerical Simulation (INS),
University of Bonn,
Wegelerstr. 6-8,
D-53211 Bonn
}
\date{}
\begin{abstract}
We show that the $d$-dimensional Haar system $H^d$ on the unit cube $I^d$ is a Schauder basis in the classical Besov space $\mathbf{B}_{p,q,1}^s(I^d)$, $0<p<1$, defined by first order differences in the limiting case $s=d(1/p-1)$,
if and only if $0<q\le p$. For $d=1$ and $p<q<\infty$, this settles the only open case in our 1979 paper \cite{Os1979}, where the Schauder basis property of $H$ in $\mathbf{B}_{p,q,1}^s(I)$ for $0<p<1$ was left undecided. 
We also consider the Schauder basis property of $H^d$ for the standard Besov spaces $B_{p,q}^s(I^d)$ defined by Fourier-analytic methods in the limiting cases $s=d(1/p-1)$ and $s=1$, complementing results by Triebel \cite{Tr1978}.
\end{abstract}
\end{frontmatter}

\section{Introduction}\label{sec1} The classical univariate Haar system $H:=\{h_m\}_{m\in\mathbb{N}}$ was one of the first examples 
of a Schauder basis in some classical function spaces on the unit interval $I:=[0,1]$. In this note, we deal with various Besov spaces ${B}_{p,q}^s$ on the unit cube $I^d\subset \mathbb{R}^d$
for the parameter range 
\be\label{PR}
0<p\le 1, \qquad 0<q<\infty,\qquad s>0,
\ee
and complement early results by Triebel \cite{Tr1978} and this author \cite{Os1979} by settling the remaining limiting cases, where the Schauder basis property of the multivariate Haar system $H^d$ was not known until now (for detailed definitions, we refer to the following sections).

There are many alternative definitions (Fourier-analytic, local means, atoms, approximations, differences, ...) that may lead to different Besov spaces for certain parts of the parameter range (\ref{PR}), see, e.g., \cite{Tr2006} for a brief introduction to function spaces of Besov-Hardy-Sobolev spaces on $\mathbb{R}^d$ and on domains.  
We consider the by now standard Besov spaces $B_{p,q}^s(I^d)$ of distributions defined in terms of Littlewood-Paley type norms (or equivalently, in terms of
atomic decompositions or local means), and the classical Besov spaces $\mathbf{B}_{p,q,1}^s(I^d)$ of functions defined by first-order differences (or, equivalently,
by best approximations with dyadic step functions). In the parameter range (\ref{PR}), these two scales of
Besov spaces coincide up to equivalent quasi-norms if and only if
\be\label{PR1}
d/(d+1)<p\le 1, \quad 0<q<\infty,\quad d(1/p-1)<s<1.
\ee

In \cite{Tr1978} it was proved that $H^d$ forms a Schauder basis in
$B_{p,q}^s(I^d)$ in the parameter range (\ref{PR1}),
see also \cite[Section 1.7.2]{Tr2006} and \cite[Section 2.5.1]{Tr2008},
where additionally the unconditionality of the Haar basis was established.
Moreover, it was also shown that the Haar system cannot be a Schauder basis
in $B_{p,q}^s(I^d)$, $0<p\le 1$, $0<q<\infty$, if either $s<d(1/p-1)$
or $s>1$.
Recently, there has been renewed interest in investigating low-order spline wavelet systems such as the Haar and Faber-Schauder systems and their multivariate counterparts as bases in Besov-Hardy-Sobolev spaces on $\mathbb{R}^d$ and $I^d$. We refer e.g. to \cite{Tr2008,Tr2010,GSU2018} and the many references cited therein. However, for $B_{p,q}^s(I^d)$ the limiting cases $s=d(1/p-1)$ and $s=1$, which were not settled in \cite{Tr1978}, are still open. We also mention the recent paper \cite{YSY2018} directly related to this note, where the authors study necessary and sufficient conditions on the parameters $p,q,s,\tau$ for which the map
$f \to (f,\chi_{I^d})_{L_2}=\int_{I^d} f\,dx$ extends to a bounded linear functional on Besov-Morrey-Campanato-type spaces $B_{p,q}^{s,\tau}(\mathbb{R}^d)$.

As to Besov spaces defined by differences,
in \cite{Os1979} it was shown that for $0<p<1$ the univariate Haar system $H$ is a Schauder basis in $\mathbf{B}_{p,q,1}^s(I)$ if
$1/p-1<s<1/p$ and $0<q<\infty$.  If $0<s<1/p-1$, then $\mathbf{B}_{p,q,1}^s(I)$ has a trivial dual space, and thus cannot possess a Schauder basis. 
For $s\ge 1/p$, $0<q<\infty$, the spaces $\mathbf{B}_{p,q,1}^s(I)$ degenerate to containing only constant functions. In the only remaining limiting case $s=1/p-1$, the proof in \cite{Os1979} established the Schauder basis property of the Haar system $H$
in $\mathbf{B}_{p,q,1}^s(I)$ also for $0 < q \le p$ while for the parameter range $p<q<\infty$ the question was left open. 

Our main goal in this paper is to settle the limiting cases for both scales of Besov spaces. In Section \ref{sec2} we will prove the following:

\begin{theo} \label{theo1} Let $d=1,2,\ldots$, and let $p,q,s$ satisfy (\ref{PR}). The Haar system $H^d$ (equipped with any of its natural enumerations) is a Schauder basis in the Besov space $\mathbf{B}_{p,q,1}^s(I^d)$ if either
\be\label{PR2}
d(1/p-1)<s<1/p, \quad (d-1)/d<p\le 1, \quad 0<q<\infty.
\ee
or if
\be\label{PR3}
s=d(1/p-1), \quad (d-1)/d<p<1, \quad 0<q\le p.
\ee
In all other cases, $H^d$ cannot be a Schauder basis in $\mathbf{B}_{p,q,1}^s(I^d)$. More precisely:\\
\bi
\item[i)] Let $0<q<\infty$. If $0<s<\min(d(1/p-1), 1/p)$ then $\mathbf{B}_{p,q,1}^s(I^d)$ has a trivial dual, while for $s\ge 1/p$ it degenerates to containing only constant functions.\\
\item[ii)] If $s=d(1/p-1)$, $(d-1)/d<p<1$, then we have two cases: \\
\bi
\item[a)] If $1<q<\infty$ then the coefficient functionals of the Haar expansion which are uniquely defined on $\mathrm{span}(H^d)$ cannot be extended to bounded linear functionals on $\mathbf{B}_{p,q,1}^{d(1/p-1)}(I^d)$. 
\item[b)] If $p<q\le 1$ then the partial sum operators of the Haar expansion
are not uniformly bounded on $\mathbf{B}_{p,q,1}^{d(1/p-1)}(I^d)$. 
\ei
\ei
\end{theo}
For $d=1$,  the statement of Theorem \ref{theo1} except for part ii)
has been established in \cite{Os1979} using characterizations of  
$\mathbf{B}_{p,q,1}^{s}(I)$ in terms of best approximations by dyadic step functions. This approach carries over to the case $d>1$.
The assertions in part ii) are new, and follow from modifying the univariate examples used in \cite{Os1979} (see the lemma on p. 535 there).

As will be clear from the proofs, the formulation of Theorem \ref{theo1} 
carries over to the Haar system on $\mathbb{R}^d$ and the Besov spaces
$\mathbf{B}_{p,q,1}^s(\mathbb{R}^d)$ without change. 
Similar results are expected to hold for Besov spaces $\mathbf{B}_{p,q,r}^{d(1/p-1)}(I^d)$ defined in terms of $r$-th order differences, $r>1$,
and multivariate spline systems of higher order, such as the Franklin system.
As to assertion ii) a), we do not know whether $\mathbf{B}_{p,q,1}^{d(1/p-1})(I^d)$ has a nontrivial dual for $q>1$ at all. 

In Section \ref{sec3}, we deal with the standard Besov spaces $B_{p,q}^{s}(I^d)$
and use their characterizations in terms of atomic decompositions and local means
to prove the following result.
\begin{theo} \label{theo2} Let $d=1,2,\ldots$, and let $p,q,s$ satisfy (\ref{PR}), where additionally $d(1/p-1)\le s \le 1$. The Haar system $H^d$ is a Schauder basis in the Besov space $B_{p,q}^s(I^d)$ if either (\ref{PR1}) or 
\be\label{PR4}
s=d(1/p-1), \quad d/(d+1) < p<1, \quad 0<q\le p.
\ee
holds.
In all other cases,
$H^d$ cannot be a Schauder basis in $B_{p,q}^s(I^d)$.  In particular,
\bi
\item[i)] If $s=1$, $d/(d+1)\le p<1$, $0<q<\infty$, then the Haar expansion of the smooth function 
$f(x)=x_1+\ldots +x_d$ does not converge to $f$ in ${B}_{p,q}^s(I^d)$\\
\item[ii)] If $s=d(1/p-1)$, $d/(d+1)\le p <1$, $p<q<\infty$, we have again two cases. \\
\bi
\item[a)] If $1<q<\infty$, then the coefficient functionals of the Haar expansion which are uniquely defined on $\mathrm{span}(H^d)$ cannot be extended to bounded linear functionals on ${B}_{p,q}^{d(1/p-1)}(I^d)$. 
\item[b)] If $p<q\le 1$, then the partial sum operators of the Haar expansion
are not uniformly bounded on ${B}_{p,q}^{d(1/p-1)}(I^d)$. 
\ei
\ei
\end{theo} 
Compared to \cite{Tr1978} only the proof of the Schauder basis property
for the parameter range (\ref{PR4}), the limiting case $s=1$ in part i), and
part ii) of Theorem \ref{theo2} are new. The theorem holds for the suitably enumerated Haar system on $\mathbb{R}^d$ and the spaces 
${B}_{p,q}^{s}(\mathbb{R}^d)$ without changes in the formulation.
The result of case a)  in part (ii) is also covered by \cite[Corollary 2.7, (ii)]{YSY2018}.

\medskip\noindent
{\bf Acknowledgments}. This research grew out of discussions
with T. Ullrich who shared the preprint version of the paper \cite{GSU2018} with me, and encouraged me to reconsider the open problems in \cite{Os1979}.  T. Ullrich also pointed out the counterexamples for part ii) in Theorem \ref{theo2} which are reproduced here with his permission. His contributions and interest are gratefully acknowledged. We also thank W. Sickel who made an early version of \cite{YSY2018} available to us.

The work on this paper started during a recent stay at the Institute for
Numerical Simulation (INS) of the University of Bonn sponsored by the Hausdorff Center for Mathematics and the Deutsche Forschungsgemeinschaft (DFG). I
thank my colleagues at the INS for the fruitful and friendly atmosphere, and the above named institutions for their support.

\section{Proof of Theorem \ref{theo1}}\label{sec2}

\subsection{Definitions and preparations}\label{sec21}
Recall first the definition of the $L_2$-normalized Haar functions. By $\chi_\Omega$ we denote the characteristic function of a Lebesgue measurable set $\Omega\subset \mathbf{R}^d$, 
and by $\Delta_{k,i}=[(i-1)2^{-k},i2^{-k})$, $i=1,\ldots,2^k$, the univariate dyadic intervals of length $2^{-k}$, $k=0,1,\ldots$. Then the univariate Haar system
$H=\{h_m\}_{m\in\mathbb{N}}$ on $I$ is given by
$h_1=\chi_I$, and
$$
h_{2^{k-1}+i}=2^{(k-1)/2}(\chi_{\Delta_{k,2i-1}}-\chi_{\Delta_{k,2i}}), \qquad i=1,\ldots,2^{k-1}, \quad k\in\mathbb{N}.
$$
As is well known, the Haar functions $h_m$ with $m\ge 2$ can also be indexed
by their supports, and identified with the appropriately scaled shifts and dilates
of a single function, the centered Haar wavelet
\be\label{HW0}
h_0:=\chi_{[-1/2,0]}-\chi_{[0,1/2]}.
\ee
Indeed,
$$
h_{\Delta_{k-1,i}}:=h_{2^{k-1}+i}=|\Delta_{k-1,i}|^{-1/2} h_0(2^{k-1} \cdot-i+1/2)
$$
for $i=1,\ldots,2^{k-1}$, and $k\in\mathbb{N}$. The above introduced enumeration of the Haar functions $h_m$ is the natural ordering used in the literature,
however, one can also define $H$ as the union of dyadic blocks
$$
H=\cup_{k=0}^\infty H_k,\qquad H_0=\{h_1\},\quad H_k=\{h_{\Delta_{k-1,i}}:\,i=1,\ldots,2^{k-1}\},\quad k\in \mathbb{N},
$$
and allow for arbitrary orderings within each block $H_k$.  
Below, we will work with the multivariate counterparts of the spaces
$$
S_k=\mathrm{span}(\{h_m\}_{m=1}^{2^k})=\mathrm{span}(\{\chi_{\Delta_{k,i}}\}_{i=1}^{2^k}),\qquad k=0,1,\ldots,
$$
of piecewise constant functions with respect to the uniform dyadic partition $T_k=\{\Delta_{k,i}:\,i=1,\ldots,2^k\}$ of step-size $2^{-k}$ on the unit interval $I$,
which we call for short dyadic step functions of level $k$.

Consider now the isotropic multivariate Haar system $H^d$ on the $d$-dimensional cube $I^d$, $d>1$, which 
we define in a blockwise fashion as follows. Let the partition $T_k^d$ be the set of all dyadic cubes of side-length $2^{-k}$ in $I^d$. Each cube in $T_k^d$ is the $d$-fold 
product of univariate $\Delta_{k,i}$, i.e., 
$$
T_k^d =\{ \Delta_{k,\mathbf{i}}:=\Delta_{k,i_1}\times \ldots\times \Delta_{k,i_d}:\;\mathbf{i}=(i_1,\ldots,i_d)\in \{1,\ldots,2^k\}^d\}.
$$
The set of all piecewise constant functions on $T_k^d$ is denoted by $S_k^d$.
With each $\Delta_{k-1,\mathbf{i}}\in T_{k-1}^d$, $\mathbf{i}\in \{1,\ldots,2^{k-1}\}^d$, $k\in \mathbb{N}$,
we can associate a set $H^d_{k,\mathbf{i}}\subset S_k^d$ of $2^d-1$ multivariate Haar functions with support $\Delta_{k-1,\mathbf{i}}$, given by all possible tensor products
$$
\psi_{k,i_1}(x_1)\cdot\psi_{k,i_2}(x_2)\cdot \ldots \cdot\psi_{k,i_d}(x_d),\qquad \psi_{k,i} = h_{\Delta_{k-1,i}}\mbox{ or } 2^{(k-1)/2}\chi_{\Delta_{k-1,i}}
$$
where at least one of the $\psi_{k,i}$ equals $h_{\Delta_{k-1,i}}$.

The blocks $H_k^d$ that define the $d$-dimensional Haar system 
$$
H^d = \cup_{k=0}^\infty H_k^d
$$
are given as follows:
The block $H_0^d$ is exceptional, and  consists of the single constant function $\chi_{I^d}$. The block $H_1^d$ coincides with 
$H^d_{1,\mathbf{1}}$ and consists of $2^d-1$ Haar functions (we use the notation $\mathbf{1}=(1,\ldots,1)$, $\mathbf{2}=(2,\ldots,2)$, etc.). For general $k\ge 2$, the block
$$
H_k^d := \cup_{\Delta_{k-1,\mathbf{i}}\in T_{k-1}^d} \, H^d_{k,\mathbf{i}}
$$
consists of $(2^d-1)2^{(k-1)d}$ Haar functions which we call Haar functions of level $k$. It is obvious that
$$
S_k^d = \mathrm{span}(\cup_{l=0}^k H_l^d),
$$
and that $H^d$ is a complete orthonormal system in $L_2(I^d)$. Since each Haar function in $H^d$ has support on a $d$-dimensional cube,
we call this system isotropic Haar system (in contrast to the $d$-dimensional tensor-product Haar system, where the supports of the 
Haar functions are $d$-dimensional dyadic rectangles).
As mentioned before for the univariate case, the ordering of the multivariate Haar functions within the blocks $H^d_k$ can be arbitrary. The statements of Theorems \ref{theo1} and \ref{theo2} hold for any enumeration of $H^d$ as long as the enumeration
does not violate the natural ordering by level $k$.

The Besov spaces $\mathbf{B}_{p,q,1}^s(I^d)$ considered in this section are defined for $s>0$, $0<p,q\le \infty$, as the set of all
(equivalence classes of) Lebesgue measurable functions $f:\,I^d\to \mathbb{R}$ for which the quasi-norm
$$
\|f\|_{\mathbf{B}_{p,q,1}^s}:=\|f\|_{L_p(I^d)}+ \| t^{-s-1/q} \omega(t,f)_p\|_{L_q(I)}
$$
is finite. Here, 
$$
\omega(t,f)_p:=\sup_{0<|y|\le t} \|\Delta_y f\|_{L_p(I^d_y)}, \qquad t>0,
$$
stands for the first-order $L_p$ modulus of smoothness, where
$$
\Delta_y f(x):=f(x+y)-f(x), \qquad x\in I^d_y:=\{z\in I^d:\;z+y\in I^d\}, \quad y\in \mathbb{R}^d,
$$ 
denotes the first-order forward difference. Here and throughout the paper, we adopt the following notational convention: If the domain is $I^d$,
we omit the domain in the quasi-norm notation, e.g., we write $\|\cdot\|_{L_p}$ instead of
$\|\cdot\|_{L_p(I^d)}$. Also, by $c,C$ we denote generic positive constants that may change from line to line, and, unless stated otherwise, depend on $p,q,s$ only. The notation $A\approx B$ is used if $cA\le B\le CA$ holds for two such constants $c,C$.

For the case $0<p<1$, $0<q<\infty$ we are interested in, $\mathbf{B}_{p,q,1}^s(I^d)$ is a  quasi-Banach space equipped with
a $\gamma$-quasi-norm, where $\gamma=\min(p,q)$, meaning that $\|\cdot\|_{\mathbf{B}_{p,q,1}^s}$
is homogeneous and satisfies
$$
\|f+g\|_{\mathbf{B}_{p,q,1}^s}^{\gamma}\le \|f\|_{\mathbf{B}_{p,q,1}^s}^{\gamma}+\|g\|_{\mathbf{B}_{p,q,1}^s}^{\gamma}.
$$
Similarly, the $L_p$ quasi-norm is a $p$-quasi-norm if $0<p\le 1$. All spaces introduced
in the sequel have $\gamma$-quasi-norms for some suitable $\gamma\in (0,1]$.

For the parameter region (\ref{PR}), the spaces are nontrivial only if $s<1/p$.
Indeed, if $f\in  \mathbf{B}_{p,q,1}^s(I^d)$ for some $s\ge 1/p$ then using the properties of the first-order $L_p$ modulus of smoothness we have $\omega(t,f)_p=\mathrm{o}(t^{1/p})$, $t\to 0$, which in turn implies $\omega(t,f)_p=0$ for all $t>0$ and $f(x)=\xi$ for some constant $\xi\in\mathbb{R}$ almost everywhere on $I^d$. From now on, we can therefore 
restrict ourselves to $0<s<1/p$ in (\ref{PR}).

In this section we will exclusively work with an equivalent quasi-norm based on approximation techniques using piecewise constant approximation on dyadic partitions.  Let
$$
E_k(f)_p:= \inf_{s\in S_k^d} \|f-s\|_{L_p},\qquad k=0,1,\ldots,
$$
denote the best approximations to $f\in L_p(I^d)$ with respect to $S_k$. From \cite[Theorem 6]{Os1980} for $d=1$, and \cite[Theorem 5.1]{DP1988} for $d>1$ we have that
\be\label{Bnorm}
\|f\|_{\mathbf{A}_{p,q,1}^s}:=\|f\|_{L_p} + (\sum_{k=0}^\infty (2^{ks} E_k(f)_p)^q)^{1/q}
\ee
provides an equivalent quasi-norm on $\mathbf{B}_{p,q,1}^s(I^d)$ for all $0<p<1$, $0< q< \infty$, $0<s<1/p$. This norm equivalence automatically implies that the set of all dyadic step functions 
\be\label{Span}
S^d:=\mathrm{span}(H^d)=\mathrm{span}(\{S_k^d\}_{k=0}^\infty) 
\ee
is dense in $\mathbf{B}_{p,q,1}^s(I^d)$ for all those parameter
values. Note that in \cite{DP1988} the case $1\le s<1/p$ is formally excluded in the formulations but the proofs in \cite{DP1988} extend to this parameter range as well. 

At the heart of the counterexamples used for the proof of Theorem \ref{theo1} is the following simple observation
which we formulate as
\begin{lem}\label{lem1}
Let $(\Omega,\mathcal{A},\mu)$ be a finite measure space, and $f\in L_p(\Omega):=L_p(\Omega,\mathcal{A},\mu)$, $0<p\le 1$, be supported on $\Omega'\in \mathcal{A}$, where  $\mu(\Omega')\le \frac12 \mu(\Omega)$. Then
$$
\|f\|_{L_p(\Omega)} = \inf_{\xi\in\mathbb{R}} \|f-\xi\|_{L_p(\Omega)},
$$
i.e., best approximation by constants in $L_p(\Omega)$ is achieved by setting $\xi=0$.
\end{lem}
{\bf Proof}. Indeed, under the above assumptions and by the inequality $|a+b|^p\le |a|^p+|b|^p$ we have
\bea
\|f-\xi\|_{L_p(\Omega)}^p&=&\int_{\Omega'}|f(x)-\xi|^p\,d\mu(x) +\mu(\Omega\backslash \Omega')|\xi|^p \ge \int_{\Omega'}(|f(x)-\xi|^p +|\xi|^p)\,d\mu(x)\\
&\ge&
\int_{\Omega'}|f(x)|^p\,d\mu(x)=\|f\|_{L_p(\Omega)}^p
\eea
for any $\xi \in \mathbb{R}$, with equality for $\xi=0$. This gives the statement. \hfill $\Box$

\medskip
Note that the equivalence (up to constants depending on parameters but not on $f$) between $L_p$ quasi-norms and 
best approximations by constants holds also for $p\ge 1$ and under weaker assumptions on the support size
of $f$ (e.g., $\mu(\Omega')/ \mu(\Omega)\le \delta<1$ would suffice). We will apply this lemma to our examples of dyadic step functions constructed below, and to the Lebesgue measure on dyadic cubes in $I^d$, where the step functions are not constant. Extensions to higher degree polynomial and spline approximation are possible as well (see the proof of the Lemma on p. 535 in \cite{Os1979} for $d=1$).

If $H^d$ is a Schauder basis in a quasi-Banach space $X$ of functions or distributions
defined on $I^d$ then necessarily $S^d=\mathrm{span}(H^d)\subset X$ and any
dyadic step function $g\in S^d$ has a unique Haar expansion 
given by
\be\label{HE}
g(x)= \sum_{h\in H^d} c_h(g) h(x),\qquad c_h(g):=\int_{I^d} g h \,dx.
\ee
Since for $g\in S^d$ only finitely many coefficients $c_h(g)$
do not vanish, the summation in (\ref{HE}) is finite, and there are no convergence issues.
Thus, for the Schauder basis property of $H^d$ in $X$ to hold, the coefficient functionals $c_h(g)$ must be extendable to elements in $X'$, and the level $k$ partial
sum operators
\be\label{Pk}
(P_kg)(x) = \sum_{l=0}^k \sum_{h\in H^d_l} c_h(g) h(x), \qquad k=0,1,\ldots,
\ee
must be extendable to \emph{uniformly} bounded linear operators in $X$. When applied to the case $X=\mathbf{B}_{p,q,1}^{s}(I^d)$ considered in this section, this explains that the statements in Theorem \ref{theo1} i)-ii) imply the failure of the Schauder basis property
of $H^d$ in $\mathbf{B}_{p,q,1}^{s}(I^d)$ for the associated parameter ranges.
The same is true for $X=B_{p,q}^{s}(I^d)$ and Theorem \ref{theo2} considered in
Section \ref{sec3}. 

For $X=L_1(I^d)$, $P_k$ extends to a bounded projection on $L_1(I^d)$ with range $S_k^d$, and with constant values on the dyadic cubes in $T_k^d$ explicitly given by averaging. This comes in handy when computing $P_kf$ for concrete functions $f$. Indeed, the constant values taken by $P_kf$ on dyadic intervals in $T_k^d$ are given by 
\be\label{AvP}
P_kf(x)=av_\Delta(f), \qquad x\in \Delta,\quad \Delta\in T_k^d,\qquad f\in L_1(I^d),
\ee
where the functionals
\be\label{Av}
av_\Delta(f):=2^{kd}\int_{I^d} \chi_\Delta f\,dx=2^{kd}\int_\Delta f\, dx,\qquad \Delta\in T_k^d,\quad k=0,1,\ldots,
\ee 
represent the average value of $f$  on dyadic cubes (we will call them for short average functionals). Note that coefficient functionals $c_h$ are finite linear combinations of average functionals as defined in (\ref{Av}), vice versa. Finally, for $X=L_2(I^d)\subset L_1(I^d)$ the level $k$ partial sum operator $P_k$ realize the orthoprojection onto $S_k^d$.

\subsection{Proof of Theorem \ref{theo1}: Positive results} \label{sec22}
For $d=1$, the cases in Theorem \ref{theo1}, where the Schauder basis property can be established, are covered by \cite{Os1979}. The proof for $d>1$  does not hold surprises,
we give it for completeness. By density of $S^d=\mathrm{span}(H^d)$ in the quasi-Banach space $\mathbf{B}_{p,q,1}^s(I^d)$, it is sufficient to establish the inequality
\be\label{Unif}
\|Pg\|_{\mathbf{A}_{p,q,1}^s}\le C\|g\|_{\mathbf{A}_{p,q,1}^s},\qquad
g\in S^d,
\ee
for any partial sum operator $P$ of the Haar expansion (\ref{HE}), with a constant $C$ independent of $g$ and $P$, for the parameters satisfying (\ref{PR2}) or (\ref{PR3}).
According to our ordering convention for $H^d$, any such partial sum operator $P$
can be written, for some $k=0,1,\ldots$ and some subset $\bar{H}^d_{k+1}\subset H^d_{k+1}$, in the form
\be\label{defPS}
{P}g=P_kg + \sum_{h\in \tilde{H}^d_{k+1}} c_h(g)h \in S_{k+1}^d.
\ee
For $\bar{H}^d_{k+1}=\emptyset$, we get $P=P_k$ as partial case.

The first step for establishing (\ref{Unif}) is the proof of the inequality
\be\label{Step11}
\|{P}g\|_{L_p}^p\le C 2^{kd(p-1)} \sum_{\Delta\in T_k^d} \|g\|_{L_1(\Delta)}^p,
\ee
with the explicit constant $C=2^d$.
By (\ref{AvP}) and (\ref{Av}), we have
$$
\|P_kg\|_{L_p}^p =\sum_{\Delta\in T_k^d} 2^{-kd} \left(2^{kd}\int_\Delta g \,dx\right)^p\le 2^{kd(p-1)} \sum_{\Delta\in T_k^d} \|g\|_{L_1(\Delta)}^p.
$$
The remaining $h\in \bar{H}^d_{k+1}$ can be grouped by their supports (which are dyadic cubes $\Delta\in T_k^d$ by construction), where each group may hold up to $2^d-1$
Haar functions with the same $\mathrm{supp}(h)=\Delta\in T_k^d$. Individually, by the definition of the Haar coefficients $c_h(g)$ and the scaling $\|h\|_{L_\infty(\Delta)}=2^{kd/2}$ of the Haar functions in $H_{k+1}^d$, we obtain for each term associated with 
a Haar function in such a group the estimate
$$
\|c_h(g)h\|_{L_p}^p = |c_h(g)|^p \|h\|_{L_p(\Delta)}^p \le 2^{kdp/2}\|g\|_{L_1(\Delta)}^p\cdot 2^{-kd}2^{kdp/2}=2^{kd(p-1)}\|g\|_{L_1(\Delta)}^p.
$$
Thus, applying the $p$-quasi-norm triangle inequality for $L_p(I^d)$ in the appropriate way, we obtain (\ref{Step11}).

Now we apply the embedding theorem $\mathbf{B}_{p,p,1}^{d(1/p-1)}(I^d)\subset L_1(I^d)$, with the appropriate coordinate transformation, to the terms $\|g\|_{L_1(\Delta)}^p$
(the stronger, optimal embedding $\mathbf{B}_{p,1,1}^{d(1/p-1)}(I^d)\subset L_1(I^d)$ is covered by \cite[Theorem 7.1]{DP1988}). This gives
$$
\|g\|_{L_1(\Delta)}^p\le C2^{kd(1-p)}(\|g\|_{L_p(\Delta)}^p+\sum_{l=0}^\infty 2^{ld(1-p)}E_{k+l}(g)_{p,\Delta}^p)
$$
for each $\Delta\in T_k^d$, where 
$$
E_{k+l}(g)_{p,\Delta}:=\inf_{s\in S_{k+l}^d} \,\|g-s\|_{L_p(\Delta)},\qquad l=0,1,\ldots,
$$
denotes the local best $L_p$ approximation by dyadic step functions restricted to cubes $\Delta$ from $T_k^d$. Since
$$
\|g\|_{L_p}^p=\sum_{\Delta\in T_k^d} \|g\|_{L_p(\Delta)}^p,\qquad 
E_{k+l}(g)_{p}^p=\sum_{\Delta\in T_k^d} E_{k+l}(g)_{p,\Delta}^p,\quad l=0,1,\ldots,
$$
after substitution into (\ref{Step11}), we arrive at the estimate
\be\label{Step21}
\|{P}g\|_{L_p}^p\le C(\|g\|_{L_p}^p+2^{kd(p-1)}\sum_{l=k}^\infty 2^{ld(1-p)}E_{l}(g)_{p}^p)
\ee
for the $L_p$ quasi-norm of any partial sum $Pg$.

With the auxiliary estimate (\ref{Step21}) at hand, we turn now to the estimate of the Besov quasi-norm of $g-Pg$. Since $Pg\in  S_{k+1}^d$, we have
$$
E_{l}(g-Pg)_p=E_{l}(g)_p,\qquad l>k,
$$
while for $l\le k$ the trivial bound 
$$
E_{l}(g-Pg)_{p}\le \|g-Pg\|_{L_p(I^d)}
$$
will suffice. This gives
$$
\|g-Pg\|_{\mathbf{A}_{p,q,1}^{s}}=\|p-Pg\|_{L_p} + (\sum_{l=0}^\infty (2^{ls} E_l(g-Pg)_{p})^q)^{1/q}\qquad\quad
$$
\be\label{Step31}
\qquad\qquad\qquad\le C\left(2^{ks}\|g-Pg\|_{L_p} + (\sum_{l=k+1}^\infty (2^{ls} E_l(g)_{p})^q)^{1/q}\right),
\ee
uniformly for all $P$ and $g\in S^d$. 

To deal with the term $\|g-Pg\|_{L_p}$, we introduce the element $s_k\in S_k^d$ of best $L_p$ approximation, i.e.,
$$
\|g-s_k\|_{L_p}=E_k(g)_{p},
$$
and estimate with (\ref{Step21}) and $Ps_k=P_ks_k=s_k$ as follows:
\bea
\|g-Pg\|_{L_p}^p&\le& \|g-s_k\|_{L_p}^p+\|P(g-s_k)\|_{L_p}^p\\
&\le& \|g-s_k\|_{L_p}^p+C(\|g-s_k\|_{L_p(I^d)}^p+2^{kd(p-1)}\sum_{l=k+1}^\infty 2^{ld(1-p)}E_{l}(g-s_k)_{p}^p)
\eea
\be\label{Step3a}
\;\;\qquad\le C2^{kd(p-1)}\sum_{l=k}^\infty 2^{ld(1-p)}E_{l}(g)_{p}^p.
\ee

Thus, since $s=d(1/p-1)$, $(d-1)/d<p<1$, and $q/p\le 1$  according to (\ref{PR3}), for the first term in the right-hand side of  (\ref{Step31}) we get
$$
2^{kd(1/p-1)}\|g-Pg\|_{L_p}\le C\left(\left(\sum_{l=k}^\infty 2^{ld(1-p)}E_{l}(g)_{p}^p
\right)^{q/p}\right)^{1/q}\le C\left(\sum_{l=k}^\infty (2^{ld(1/p-1)}E_{l}(f)_{p})^q\right)^{1/q},
$$
where the inequality
\be\label{pinq}
\sum_{l=0}^\infty a_l \le (\sum_{l=0}^\infty a_l^\gamma)^{1/\gamma}, \qquad a_l\ge 0, \quad 0<\gamma \le 1,
\ee
has been used with $\gamma=q/p$, $a_l=2^{ld(1-p)}E_{l}(g)_{p}^p$ for $l\ge k$, and $a_l=0$ for $l<k$.
After substitution into
(\ref{Step3}) we arrive at
\be\label{Unif1}
\|g-Pg\|_{\mathbf{A}_{p,q,1}^{d(1/p-1)}}\le C\left(\sum_{l=k}^\infty (2^{ls} E_l(g)_{p})^q\right)^{1/q}\le C\|g\|_{\mathbf{A}_{p,q,1}^{d(1/p-1)}}
\ee
for all $g\in S^d$ if the parameters satisfy (\ref{PR3}).
Since the quasi-norm in $\mathbf{A}_{p,q,1}^{s}(I^d)$ is a $\min(p,q)$-quasi-norm,
(\ref{Unif1}) is equivalent with (\ref{Unif}). This proves the Schauder basis
property for $H^d$ in $\mathbf{B}_{p,q,1}^{s}(I^d)$ for the parameters satisfying (\ref{PR3}).

For the parameter range (\ref{PR2}), i.e., when $d(1/p-1)<s<1/p$, $(d-1)/d<p\le 1$, $0<q<\infty$, we can  apply the Hardy-type inequality 
$$
(\sum_{l=k}^\infty a_l^p)^{1/p}\le C_{\epsilon,q/p}2^{-k\epsilon}(\sum_{l=k}^\infty (2^{l\epsilon}a_l)^q))^{1/q},\quad \epsilon >0, \quad k=0,1,\ldots,
$$
valid for non-negative sequences $\{a_l\}$ and all $0 <p,q <\infty$. Setting $\epsilon=s-d(1/p-1)$ and $a_l=2^{ld(1/p-1)}E_{l}(f)_{p}$,
from (\ref{Step3a}) we obtain
$$
2^{ks}\|g-Pg\|_{L_p}\le C 2^{k(s-d(1/p-1))}(\sum_{l=k}^\infty 2^{ld(1-p)}E_{l}(g)_{p}^p)^{1/p}\le C(\sum_{l=k}^\infty 2^{lds}E_{l}(g)_{p}^q)^{1/q}.
$$
It remains to substitute this into (\ref{Step31}) and proceed as above for the parameter
range (\ref{PR3}).
This concludes the proof of the Schauder basis property for all parameters satisfying (\ref{PR2}) or (\ref{PR3}).

\subsection{Proof of Theorem \ref{theo1}: Negative results}\label{sec23}
 We first deal with part i),
and follow the proof given in \cite{Os1979} for $d=1$. The case $s\ge 1/p$ has been discussed before. Let $0<s<d(1/p-1)$, $0<p<1$, $0<q<\infty$,
and assume that $F$ is a bounded linear functional on $\mathbf{B}_{p,q,1}^{s}(I^d)$. For any dyadic cube $\Delta\in T_k^d$
of side-length $2^{-k}$ we have by Lemma \ref{lem1}
$$
\|\chi_\Delta\|_{L_p}=E_l(\chi_\Delta)_p=2^{-kd/p},\quad l=0,\ldots,k-1,\qquad E_l(\chi_\Delta)_p=0,\quad l=k,k+1,\ldots,
$$
and consequently
$$
\|\chi_\Delta\|_{\mathbf{A}_{p,q,1}^s}=2^{-kd/p}(1+(\sum_{l=0}^{k-1} 2^{lsq})^{1/q}
\le C 2^{k(s-d/p)},\qquad \Delta\in T^d_k,\quad k=0,1,\ldots .
$$
By linearity and boundedness of $F$ this implies
\bea
|F(\chi_\Delta)|&=& |\sum_{\Delta'\in T^d_l:\,\Delta'\subset \Delta} F(\chi_{\Delta'}|\le
\sum_{\Delta'\in T^d_l:\,\Delta'\subset \Delta} |F(\chi_{\Delta'}|\\
&\le& C \sum_{\Delta'\in T^d_l:\,\Delta'\subset \Delta} \|\chi_{\Delta'}\|_{\mathbf{A}_{p,q,1}^s}\le C2^{(l-k)d} 2^{l(s-d/p)} =C2^{-kd} 2^{l(s-d(1/p-1))},\qquad l>k,
\eea
for any given $\Delta\in T_k^d$. Here $C$ also depends on $F$. If we let $l\to \infty$, we obtain $F(\chi_\Delta)=0$
for all dyadic cubes, and by linearity $F(g)=0$ for all $g\in S^d=\mathrm{span}(H^d)$.
Since the dyadic step functions are dense in $\mathbf{B}_{p,q,1}^{s}(I^d)$,
this shows $F=0$, i.e., the dual of $\mathbf{B}_{p,q,1}^{s}(I^d)$ is trivial.

The counterexamples proving the statement in part ii) of Theorem \ref{theo1}
are new even for $d=1$ (actually, subcase b) disproves our conjecture concerning 
the parameter range $p<q\le 1$ made in \cite{Os1979}). Consider first a),
i.e., assume that $s=d(1/p-1)$, $(d-1)/d < p < 1$, and $1<q<\infty$.
Proving that the coefficient functionals in (\ref{HE}) cannot be extended boundedly
from $S^d$ to $\mathbf{B}_{p,q,1}^{{d(1/p-1)}}(I^d)$ is the same as proving this
for the average 
functionals $av_\Delta$ defined in (\ref{Av}) for arbitrary dyadic cubes and
$g\in S^d\subset L_1(I^d)$.
Since for the above parameter range the Besov space $B^{d(1/p-1)}_{p,q,1}(I^d)$ is not embedded into $L^1(I^d)$, see \cite{Os1979,Os1980,DP1988} for the corresponding embedding theorems, we see the problem.

We provide the argument for the average functional $av_{I^d}$ associated with the dyadic cube $\Delta=I^d\in T_0^d$, by showing that
there is a sequence $g_k\in S^d$ of dyadic step functions which is uniformly bounded 
in $\mathbf{B}_{p,q,1}^{{d(1/p-1)}}(I^d)$, and such that
\be\label{L1div}
av_{I^d}(g_k)=\int_{I^d} g_k \,dx \to \infty,\qquad k\to \infty.
\ee 
By dilating and shifting these $g_k$ to fit their supports into an arbitrarily given dyadic
cube $\Delta$, similar examples can be obtained for all average functionals in (\ref{Av}).
Since we do not attempt to show quantitative lower bounds for
the divergence in (\ref{L1div}), the following construction suffices for $d>1$,
the modification for $d=1$ is stated below. Set
\be\label{gka}
g_{k}= \sum_{j=1}^k a_j \chi_{\Delta_{j,\mathbf{2}}}, \qquad a_j=2^{jd} j^{-1}, \quad j\in \mathbb{N}.
\ee
This is a dyadic step function which takes positive values $a_j$ on a sequence of non-overlapping dyadic cubes $\Delta_{j,\mathbf{2}}$, $j=1,\ldots,k$, located along the main diagonal of $I^d$, and is otherwise zero. Thus,
$$
\int_{I^d} g_k \,dx =\sum_{j=1}^k 2^{-jd}a_j =\sum_{j=1}^k j^{-1}\ge c\log(k+1),
$$ 
and (\ref{L1div}) is established. 

For $d>1$ these $g_k$ also satisfy the assumptions
of Lemma \ref{lem1} with respect to any dyadic cube $\Delta$, where $g_k$ is not constant (if $g_k$ is constant on a dyadic cube, its best $L_p$ approximation
by constants on this cube is obviously zero). This allows us to compute the best approximations of $g_k$ exactly:
\bea
E_l(g_k)_p &=&\|\sum_{j=l+1}^k a_j \chi_{\Delta_{j,\mathbf{2}}}\|_{L_p} = (\sum_{j=l+1}^k a_j^p2^{-jd})^{1/p} \\
&& \left\{\ba{ll} \le C 2^{-ld(1/p-1)}(l+1)^{-1},& l=0,\ldots,k-1\\
=0,& l=k,k+1,\ldots \ea   \right.
\eea
Substituted into the expression for the $\mathbf{A}_{p,q,1}^{d(1/p-1)}$ quasi-norm, this gives
$$
\|g_k\|_{\mathbf{A}_{p,q,1}^{d(1/p-1)}}^q \le C (\sum_{l=0}^{k-1} (l+1)^{-q})^{1/q},\qquad
k=1,2,\ldots,
$$
which shows the uniform boundedness of the sequence $g_k$ in $\mathbf{A}_{p,q,1}^{d(1/p-1)}$ since $q>1$ in case a). Here, we have silently used that $\|g_k\|_{L_p}=E_0(g_k)_p$ by Lemma \ref{lem1}. For $d=1$, to enable the application of Lemma \ref{lem1} also in this case, a modified definition of the $g_k$, e.g., 
$$
g_k= \sum_{j=1}^k a_{j+1} \chi_{\Delta_{j+1,3}},
$$
will do, the details are left to the reader. Note that the above sequences $g_k$ converge to a limit function $f\in \mathbf{B}_{p,q,1}^{d(1/p-1)}(I^d)$ which does not belong to $L_1(I^d)$, for $d=1$ similar examples were used in \cite{Os1979}.

In case b), i.e., when $s=d(1/p-1)$, $(d-1)/d<p<1$, $p<q\le 1$, the coefficient functionals in (\ref{HE}) and the dyadic averaging
functionals (\ref{Av}) can be extended to bounded linear functionals
on $\mathbf{B}_{p,q,1}^{d(1/p-1)}(I^d)$. Thus, the level $k$ partial sum operators $P_k$
defined in (\ref{Pk}) can be extended to bounded operators acting in $\mathbf{B}_{p,q,1}^{d(1/p-1)}(I^d)$. However, they are not uniformly bounded as will be shown by 
a different type of examples. Fix $k=1,2,\ldots$, and consider the $2^{(k-1)d}$
dyadic cubes $\Delta_{k,\mathbf{i}}$ in $T_k^d$ 
for which all entries of $\mathbf{i}$ are odd. Select dyadic subcubes $\tilde{\Delta}_{k+j}$ of shrinking side-length $2^{-k-j}$, $j=1,2,\ldots,2^{(k-1)d}$, one in each of them. Then we define
$$
g_k=\sum_{j=1}^{2^{(k-1)d}} b_{k,j}\chi_{\tilde{\Delta}_k},\qquad
b_{k,j}=2^{(k+j)d}j^{-1/p},\quad j=1,\ldots,2^{(k-1)d}.
$$
will do. The construction of this $g_k$ is such that Lemma \ref{lem1} is again applicable, locally on each dyadic cube where $g_k$ is not constant. This allows us to compute the best approximations $E_l(g_k)_p$ as follows: For $l=0,\ldots,k$,
we have
$$
E_l(g_k)^p_p=\|g_k\|^p_{L_p}=\sum_{j=1}^{2^{(k-1)d}} 2^{-(j+k)^d}b_{k,j}^p
\le C2^{-kd(1-p)}=C2^{-ksp}.
$$
For $l=k+1,\ldots,k+2^{(k-1)d}-1$, we get similarly
$$
E_l(g_k)^p_p=\sum_{j=1+l-k}^{2^{(k-1)d}} 2^{-(k+j)^d}b_{k,j}^p 
\le C2^{-ld(1-p)}(l-k)^{-p}=C2^{-ksp}\cdot 2^{-(l-k)sp}(l-k)^{-p}(l-k)^{-1},
$$
while $E_l(g_k)^p_p=0$ for $l\ge k+2^{(k-1)d}$. 
Thus, with these formulas for the best approximations $E_l(g_k)_p$ and the substitution $j=l-k$ for $l>k$, one arrives at
\bea
\|g_k\|_{\mathbf{A}_{p,q,1}^{d(1/p-1)}}&\le& C2^{-ks}\left(1+(\sum_{l=0}^k 2^{lsq} +
2^{ksq}\sum_{j=1}^{2^{(k-1)d}-1} 2^{jsq}\cdot 2^{-jsq}j^{-q/p})^{1/q}\right)\\
&\le& C(1+(\sum_{j=1}^\infty j^{-q/p})^{1/q})\le C<\infty, \qquad k=1,2,\ldots,
\eea
since $q>p$.

On the other hand, by  (\ref{AvP}) the level $k$ partial sum $P_kg_k$ of $g_k$ is constant on dyadic cubes in $T_k^d$, and equals $2^{kd}\cdot 2^{-(k+j)d}b_{k,j}=2^{kd}j^{-1/p}$
on the cube $\Delta_{k,\mathbf{i}}$ containing $\tilde{\Delta}_{k+j}$,
$j=1,\ldots,2^{(k-1)d}$, and vanishes on all cubes $\Delta_{k,\mathbf{i}}$ for which
at least one entry in $\mathbf{i}$ is even. The latter property ensures that Lemma \ref{lem1} is also applicable to $P_kg_k$, and gives
$$
E_l(P_kg_k)^p_p=\|P_kg_k\|^p_{L_p}= \sum_{j=1}^{2^{(k-1)d}} 2^{-kd}\cdot 2^{kdp} j^{-1}
=2^{-ksp}\sum_{j=1}^{2^{(k-1)d}} j^{-1}\ge c2^{-ksp} k
$$
for $l=0,1,\ldots,k-1$. Consequently,
$$
\|g_k\|_{\mathbf{A}_{p,q,1}^{d(1/p-1)}}\ge c2^{-ks} k^{1/p}(\sum_{l=0}^{k-1} 2^{lsq})^{1/q}\ge c k^{1/p},
$$
which shows that the partial sum operators $P_k$ are not uniformly bounded on $\mathbf{B}_{p,q,1}^{d(1/p-1)}(I^d)$ for $(d-1)/d<p<q$. This concludes the proof of Theorem \ref{theo1}.

\medskip
We have not made any attempt to obtain the exact growth of norms of partial sum operators in part b) of Theorem \ref{theo1}.
For $s=d(1/p-1)$ and $(d-1)/d< p<q \le 1$, the above considerations give the lower bound
$$
\|P_k\|_{\mathbf{B}_{p,q,1}^{d(1/p-1)}\to \mathbf{B}_{p,q,1}^{d(1/p-1)}}\ge ck^{1/p},\qquad j=1,2,\ldots,
$$
which is certainly not optimal. If one takes  $b_{k,j}=2^{(k+j)d}j^{-\alpha}$ with $\alpha=1/q+\epsilon$ and small enough $\epsilon>0$ in the above definition of $g_k$ then the better estimate 
$$
\|P_k\|_{\mathbf{B}_{p,q,1}^{d(1/p-1)}\to \mathbf{B}_{p,q,1}^{d(1/p-1)}}\ge c 2^{k(1/q-1/p-\epsilon)},\qquad k=1,2,\ldots,
$$
results, where $c>0$ depends also on $\epsilon$.

\medskip
On a final note: In the literature (with the exception of \cite{Os1979,Os1980}), for $0<p<1$ the Besov spaces
$\mathbf{B}_{p,q,r}^s(I^d)$ defined as subspaces of $L_p(I^d)$ using $r$-th order moduli of smoothness are only considered 
for the values $0<s<r$ (see, e.g.,  \cite{DP1988,Tr2006}). The reason is two-fold: It is known that
for $r\le s< r+1/p-1$, $0<q<\infty$, the spaces $\mathbf{B}_{p,q,r}^s(I^d)$ are strange: With the exception of polynomials of degree $<r$, smooth functions from $C^r(I^d$) cannot belong to $\mathbf{B}_{p,q,r}^s(I^d)$, while $C^{r-2}$-smooth dyadic splines of degree $r-1$
are dense in these spaces. This is counter-intuitive, and makes their usefulness in  applications doubtful. Moreover, the spaces
${B}_{p,q}^s(I^d)$ defined by the Fourier-analytic approach that dominate the scene coincide
with $\mathbf{B}_{p,q,r}^s(I^d)$ (in the sense of having equivalent quasi-norms) only in the range $d(1/p-1)<s<r$.
In other words, our new results on the properties of the Haar system $H^d$ in Besov spaces $\mathbf{B}_{p,q,1}^s(I^d)$ for the limiting case $s=d(1/p-1)$
do not automatically answer the same question for the scale ${B}_{p,q}^s(I^d)$. The latter will be considered in the next section.

\section{Proof of Theorem \ref{theo2}}\label{sec3} 
\subsection{Definitions and preparations}\label{sec31}
The role of the Haar system as Schauder basis in the Besov spaces ${B}_{p,q}^s(I^d)$ with $0<p\le 1$, $0<q<\infty$, $s\in \mathbb{R}$, defined in Fourier-analytic terms has been examined by Triebel \cite{Tr1978} (see also \cite[Theorem 1.58]{Tr2006}) 
who settled all but the limiting cases $s=d(1/p-1)$, $d/(d+1)<p\le 1$, $0<q<\infty$, and $s=1$, $d/(d+1) \le p \le 1$, $0<q<\infty$. Theorem \ref{theo2} gives now answers
in the limiting cases as well.

The definition of ${B}_{p,q}^s(I^d)$ is reduced by restriction to the definition  of ${B}_{p,q}^s(\mathbb{R}^d)$:
\be\label{BId}
{B}_{p,q}^s(I^d)=\{f=\tilde{f}|_{I^d}: \;\tilde{f}\in {B}_{p,q}^s(\mathbb{R}^d)\},\qquad
\|f\|_{{B}_{p,q}^s}   =\inf_{\tilde{f}:\,f=\tilde{f}|_{I^d}} \|\tilde{f}\|_{{B}_{p,q}^s(\mathbb{R}^d)}.
\ee
The definition of ${B}_{p,q}^s(\mathbb{R}^d)$ will be given in terms of atoms, for the equivalent definition in Fourier-analytic terms and a short review of the various definitions for spaces on $\mathbb{R}^d$ and on domains we refer to \cite[Chapter 1]{Tr2006}. Since
we are only interested in the limiting cases of low smoothness $s=d(1/p-1)<1$
and $s=1$ in (\ref{PR}), some simplifications are possible. 
Let us go to the details. For the parameter range $d(1/p-1)\le s\le 1$ of interest, take any $\sigma>s$ (note that for $s<1$ it is always possible to take $s<\sigma\le 1$), and consider 
the set of all H{\"o}lder class functions
$a\in \mathcal{C}^\sigma(\mathbb{R}^d)$ with support in a fixed cube of side-length $C_0>1$ centered at the origin, and with $\mathcal{C}^\sigma(\mathbb{R}^d)$ norm bounded by $C_0$. Denote this set for short by $\mathcal{C}^\sigma_{C_0}$. Functions of the form
\be\label{Atoms}
a_{j,\mathbf{i}}(x)=a(2^jx-\mathbf{i}),\qquad a\in \mathcal{C}^\sigma_{C_0},\qquad \mathbf{i}\in \mathbb{Z}^d,
\ee
are called atoms of level $0$ if $j=0$, and atoms of level $j=1,2,\ldots$ if additionally
$$
\int_{\mathbb{R}^d} a_{j,\mathbf{i}}\,dx = \int_{\mathbb{R}^d} a_{j,\mathbf{i}}\,dx=0, \qquad j\in \mathbb{N}.
$$
This latter additional condition is necessary for the following statement only if $s=d(1/p-1)$, the case we are most interested in.

\begin{lem}\label{lem2} Let $d(1/p-1)\le s<\sigma $, $d/(d+1) < p < 1$, $0<q<\infty$, and $c_0>1$ be fixed. Then $\tf\in B_{p,q}^s(\mathbb{R}^d)$
if and only if
\be\label{AtRep}
\tf(x) = \sum_{j=0}^\infty \sum_{\mathbf{i}\in \mathbb{Z}^d} c_{j,\mathbf{i}}a_{j,\mathbf{i}}(x)
\ee
(unconditional convergence in $S'(\mathbb{R}^d)$) for some atoms $a_{j,\mathbf{i}}$
specified by (\ref{Atoms}) and with coefficients such that
$$
\sum_{j=0}^\infty 2^{j(s-d/p)q}\left(\sum_{\mathbf{i}\in \mathbb{Z}^d} |c_{j,\mathbf{i}}|^p\right)^{q/p}<\infty.
$$
Moreover, 
\be\label{AtNorm}
\|\tf\|^+_{B_{p,q}^s(\mathbb{R}^d)}
:=\inf\,\left(\sum_{j=0}^\infty 2^{j(s-d/p)q}
(\sum_{\mathbf{i}\in \mathbb{Z}^d} |c_{j,\mathbf{i}}|^p)^{q/p}\right)^{1/q}
\approx \|f\|_{B_{p,q}^s(\mathbb{R}^d)},
\ee
where the infimum is taken with respect to all possible representations (\ref{AtRep}),
is an equivalent quasi-norm on $B_{p,q}^s(\mathbb{R}^d)$. The constants in the norm equivalence depend on $\sigma, C_0$, and $p, q, s$.
\end{lem}

This statement is covered by \cite[Corollary 1.23 (i)]{Tr2006}, where references 
to the history of atomic characterizations of function spaces can be found. 
Note that our atoms correspond to the $1_\sigma$-atoms ($j=0$) and $(s,p)_{\sigma,1}$-atoms ($j=1,2,\ldots$) of Definition 1.21 in \cite{Tr2006} but are scaled differently.
Instead, the necessary scaling has been incorporated in the definition of the
atomic quasi-norm (\ref{AtNorm}). Below, we will apply this lemma  with values
$\sigma > 1$, and appropriately fixed $C_0$, to obtain upper bounds for
(atomic) Besov norms.

In some cases, especially for obtaining lower bounds for ${B}_{p,q}^s(\mathbb{R}^d)$ quasi-norms, it is more convenient to use characterizations in terms of local means
\be\label{LM}
\kappa(t,\tf)(x)=(\kappa^t \ast \tf)(x):=\int_{\mathbb{R}^d} \kappa^t(x-y)\tf(y)\, dy,\qquad \kappa^t(x):=t^{-d}\kappa(t^{-1}x),\quad t>0,
\ee 
where the kernel $\kappa\in C^\infty(\mathbb{R}^d)$ has support in the cube
$[-1/2,1/2]^d$, and satisfies 
\be\label{Moment}
\kappa^{\vee}(\xi)\neq 0,\quad 0<|\xi|<\epsilon,\qquad (D^\alpha\kappa^{\vee})(0)=0\;\mbox{ if }\;
|\alpha|\le s,
\ee
for some $\epsilon>0$. Here, $\kappa^{\vee}$ denotes the Fourier transform of $\kappa$.
For $s<1$, the moment condition in (\ref{Moment}) reduces to requiring
$$
\int_{\mathbb{R}^d} \kappa(x)\,dx=0,
$$
while for $s=1$ we additionally need also orthogonality to linear polynomials:
$$
\int_{\mathbb{R}^d} x_i\kappa(x)\,dx=0,\qquad i=1,\ldots,d.
$$
We also fix another kernel $\kappa_0 \in C^\infty(\mathbb{R}^d)$
with support in the cube
$[-1/2,1/2]^d$ which, instead of (\ref{Moment}), satisfies 
$$
\kappa^{\vee}_0(0)=\int_{\mathbb{R}^d} \kappa_0(x)\,dx >0.
$$ 

By applying \cite[Theorem 1.10]{Tr2006}, we have the following characterization 
in the range of parameters of interest to us.
\begin{lem}\label{lem3} Let $0 < s\le 1$, $0<p\le 1$,  $0<q<\infty$, and let the kernels $\kappa,\kappa_0$ satisfy the above conditions. Then
$$
\|\tf\|^\ast_{B_{p,q}^s(\mathbf{R}^d)}:= \left(\|\kappa_0(1,\tf)\|_{L_p(\mathbb{R}^d)}^q+\sum_{j=1}^{\infty} 2^{jsq}\|\kappa(2^{-j},\tf)\|_{L_p(\mathbb{R}^d)}^q\right)^{1/q}
\approx \|\tf\|_{ B_{p,q}^s(\mathbf{R}^d)}.
$$
The constants in the norm equivalence depend on $\kappa, \kappa_0$, and $p, q, s$.
\end{lem}

We conclude this subsection by a technical result which shows how to reduce estimates for partial sum operators $P$ associated with the Haar expansion (\ref{HE}) of functions on $I^d$ to estimates for similar operators acting on functions defined on $\mathbb{R}^d$. To this end,
to any $P$ given by (\ref{defPS}) we associate its extension 
\be\label{defTPS}
(\tP \tf)(x)=(\tP_k \tf)(x)+\sum_{h\in \bar{H}^d_{k+1}} c_h(f)h =\left\{\ba{ll} (P(\tf|_{I^d}))(x),& x\in I^d,\\ av_\Delta(\tf),& x\in \Delta\not\subset I^d,  \ea\right.
\ee
for $\tf\in L_{1,loc}(\mathbb{R}^d)$, where $\Delta$ runs through all dyadic cubes $\Delta$ of side-length $2^{-k}$ in $\mathbb{R}^d\backslash I^d$. In other words, we define $\tP$ outside $I^d$ by $\tP_k$, the natural extension of the level $k$ partial sum operator $P_k$ to functions on $\mathbb{R}^d$. Other extensions are possible, this one simplifies some considerations below.
In particular, $\tP\tf$ has the following properties which we use throughout the rest of this subsection. First of all, it is piecewise constant on dyadic cubes $\Delta\in \tT_k^d$ outside $I^d$ and $\Delta\in T_{k+1}^d$ inside $I^d$.
Here, and in the following, $\tT_k^d$ denotes the collection of all dyadic cubes
of side-length $2^{-k}$ in $\mathbb{R}^d$ (thus, $T_k^d=\tT_k^d\cap I^d$).
Moreover, in analogy to (\ref{AvP}), we have
\be\label{AvTP}
\tP_k\tf(x)=av_\Delta(\tf);\qquad x\in \Delta\in \tT_k^d,
\ee
while by the definition of the Haar functions for each $h\in H_{k+1}^d$ we have
\be\label{AvH}
|c_h(\tf)h|=\left|\sum_{\Delta^\ast\in T_{k+1}^d:\,\Delta^\ast\subset \mathrm{supp}(h)} \alpha_{h,\Delta^\ast} av_{\Delta^\ast}(\tf)\right|,\qquad \sum_{\Delta^\ast\in T_{k+1}^d:\,\Delta^\ast\subset \mathrm{supp}(h)} \alpha_{h,\Delta^\ast} =0.
\ee
Because of shift-invariance, depending on the type of Haar function $h$ and the location of the cube $\Delta^\ast$ relative to the cube $\mathrm{supp}(h)$ containing it, there appear only finitely many different coefficient sets $\{\alpha_{h,\Delta^\ast}\}$  in (\ref{AvH}). Consequently, the restriction of the difference 
$$
|(\tP-\tP_k)\tf|\le \sum_{h\in \bar{H}_{k+1}^d} |c_h(\tf)h|
$$ 
to any $\Delta^\ast\in T_{k+1}^d$ can be bounded by the sum of \emph{differences} $|av_{\Delta'}(\tf)-av_{\Delta"}(\tf)|$ of averages with respect to neighboring dyadic cubes $\Delta',\Delta''\in T_{k+1}^d$ belonging to the same dyadic cube of side-length $2^{k}$ as $\Delta^\ast$. This will be used in subsection \ref{sec33} .

From now on,  the notation $\tf\in Y$ is reserved for functions in $B_{p,q}^s(\mathbb{R}^d)$ such that $\tf\in S'(\mathbf{R}^d)$ is represneted by an atomic decomposition (\ref{AtRep}) satisfying
\be\label{CS}
c_{j,\mathbf{i}}=0 \quad \mbox{ if}\quad \mathrm{supp}(a_{j,\mathbf{i}})\subset
\mathbb{R}^d\backslash I^d.
\ee 
For a given $f$ belonging to a Besov space $B_{p,q}^s(I^d)$ for which Lemma \ref{lem2} holds, we write $\tf\in Y_f$ if $\tf\in Y$ and $\tf|_{I^d}=f$.
Then, by the definition of atomic quasi-norms we have
\be\label{Y}
\|f\|_{B_{p,q}^s}\ge c\inf_{\tilde{f}:\,f=\tilde{f}|_{I^d}} \|\tilde{f}\|^+_{{B}_{p,q}^s(\mathbb{R}^d)}=c\inf_{\tilde{f}\in Y_f\cap B_{p,q}^s(\mathbb{R}^d)} \|\tilde{f}\|^+_{{B}_{p,q}^s(\mathbb{R}^d)}. 
\ee

\begin{lem}\label{lem4} Assume that the parameters $p,q,s$ are such that Lemma \ref{lem2} and \ref{lem3} hold. Then the operator $P$ defined in (\ref{defPS}) satisfies
\be\label{BP}
\|Pf\|_{B_{p,q}^s} \le C\|f\|_{B_{p,q}^s},\qquad f\in B_{p,q}^s(I^d),
\ee
with a constant independent of $P$, if its extension $\tP$ defined in (\ref{defTPS}) satisfies
\be\label{BTP}
\|\tP\tf\|^\ast_{B_{p,q}^s(\mathbb{R}^d)} \le C\|\tf\|^+_{B_{p,q}^s(\mathbb{R}^d)},
\qquad \tf\in Y\cap B_{p,q}^s(\mathbb{R}^d),
\ee
with a constant independent of $P$.
\end{lem}

{\bf Proof.} This follows by the locality properties of partial sum operators.
By (\ref{BTP}), $\tP\tf\in B_{p,q}^s(\mathbb{R}^d)$ is meaningfully defined for $\tf\in Y_f\cap B_{p,q}^s(\mathbb{R}^d)$. Since $\tf|_{I^d}=f$, by definition of $\tP$ we also have
$$
\tP\tf|_{I^d} = Pf.
$$
Thus, $\tP\tf$ is an extension of $Pf$, and by the definition of the
$B_{p,q}^s$ quasi-norm, by Lemma \ref{lem3}, and by (\ref{BTP}) we get 
$$
\|Pf\|_{B_{p,q}^s}\le \|\tP\tf\|_{B_{p,q}^s(\mathbb{R}^d)}\le C
\|\tP\tf\|^\ast_{B_{p,q}^s(\mathbb{R}^d)}\le
C\|\tf\|^+_{B_{p,q}^s(\mathbb{R}^d)}.
$$
It remains to take the infimum with respect to $\tf\in Y_f\cap B_{p,q}^s(\mathbb{R}^d)$, and to apply (\ref{Y}). Lemma \ref{lem4} is proved.
\hfill $\Box$

\subsection{The limiting case $s=1$}\label{sec32}
We first deal with the case $s=1$, $d/(d+1) \le p\le 1$, $0<q<\infty$, 
and show that for $f(x)=x_1+\ldots +x_d$ we have
\be\label{S1}
\|f-P_kf\|_{B_{p,q}^1} \ge c >0, \qquad k=0,1,\ldots,
\ee
for some positive constant $c$. Since $f\not\in \mathrm{span}(H^d)$, we have $f-P_kf\neq 0$, and it suffices
to consider large enough $k$.

To obtain the lower bounds needed for (\ref{S1}), we compute lower estimates for the 
$L_p$ quasi-norm of $\kappa(2^{-(k+1)},\tG_k)(x)$ for any extension $\tG_k$ of $f-P_kf$ with a kernel $\kappa$ as defined in subsection \ref{sec31} for $s=1$ (in particular, $\kappa$
is orthogonal to linear polynomials (\ref{Moment})). To this end, we observe that inside $I^d$ the difference $\tG_k(x)=(f-P_kf)(x)$ coincides a.e. with the restriction to $I^d$ of a suitably dilated and scaled single integer-shift invariant function $f_0(x)$ given by
$$
f_0(x+\mathbf{i})=f(x)-d/2, \qquad x\in (0,1]^d, \quad\mathbf{i}\in\mathbb{Z^d}.
$$
Indeed, we have
$$
(f-P_kf)(x)=2^{-k}f_0(2^kx),\qquad x\in I^d,
$$
which can be checked from the formula
$$
(f-P_kf)(x)=f(x)-av_\Delta(f)=f(x)-f(x_\Delta)=2^{-k}f(2^k(x-x_\Delta)),\qquad x\in \Delta,
$$
where $\Delta$ is an arbitrary cube in $T_k^d$, and $x_\Delta$ denotes the center of $\Delta$.

Using the invariance of $f-P_kf$ with respect to shifts of the form $2^{-k}\mathbf{j}$
inside $I^d$, and the fact that $\kappa^{2^{-(k+1)}}$ has support in a cube
of side-length $2^{-(k+1)}$ centered at the origin, we see that
$$
\kappa(2^{-(k+1)},\tG_k)(x)=2^{-k}\kappa(2^{-(k+1)},f_0(2^k\cdot))(x)=2^{-k}\kappa(1/2,f_0)(2^kx)
$$
holds for all $x\in [2^{-(k+1)},1-2^{-(k+1)}]^d$.
Thus, for $k>1$ we obtain
\bea
\|\kappa(2^{-(k+1)},\tG_k)\|_{L_p(\mathbb{R}^d)}&\ge&
2^{-k} \|\kappa(1/2,f_0)(2^k\cdot)\|_{L_p([2^{-(k+1)},1-2^{-(k+1)}]^d)}\\
&\ge& 2^{-k} ((2^{k}-2)/2^k)^{d/p}\|\kappa(1/2,f_0)\|_{L_p}\ge c 2^{-k},
\eea
for some constant $c>0$ depending on the kernel and $p$. By Lemma \ref{lem3} we conclude that
$$
\|f-P_kf\|_{B_{p,q}^1}\ge c 2^{k+1} \inf_{\tG_k} \|\kappa(2^{-(k+1)},\tG_k)\|_{L_p(\mathbb{R}^d)}\ge c,\qquad k>1,
$$
which proves (\ref{S1}).

\medskip
We finish the consideration for $s=1$ with a remark concerning the special
case $p=d/(d+1)$. In the proof of (\ref{S1}) we have not made explicit use of the restriction $d/(d+1)\le p\le 1$. That (\ref{S1}) contradicts the Schauder basis property of $H^d$ in $B_{p,q}^1(I^d)$ as claimed in Theorem \ref{theo2} is clear if $d/(d+1)<p\le 1$
since for this parameter range $B_{p,q}^1(I^d)$, $0<q<\infty$, is continuously embedded into $L_1(I^d)$ which ensures that the $P_k$ are the right candidate
to be considered for partial sum operators. Moreover, the embedding also 
implies together with (\ref{S1}) that the set $S^d=\mathrm{span}(H^d)$ of dyadic step functions cannot be dense in $B_{p,q}^1(I^d)$, thus extending the similar statement for $1<s<1/p$ proved 
in \cite{Tr1978} to the case $s=1$.

For $p=d/(d+1)$, we have
$d(1/p-1)=1=s$, and the continuous embedding
$B_{d/(d+1),q}^1(I^d)\subset L_1(I^d)$ holds only if $q\le 1$.
For $q>1$, we cannot automatically exclude the possibility that there are Haar series other than (\ref{HE}) representing the above
$f$ in $B_{d/(d+1),q}^1(I^d)$. Nor do we know for sure if  $S^d$
is dense in $B_{d/(d+1),q}^1(I^d)$. However, even in this special case $H^d$
cannot be a Schauder basis since case a) in part ii) of Theorem \ref{theo2}
applies (for the proof, see the next subsection).

We finally note that an example similar to $f(x)=x_1+\ldots+x_d$ has been used in \cite[Section 4]{GSU2018} for 
showing lower bounds for $B_{p,q}^s(\mathbb{R}^d)$ quasi-norm of the level $k$ partial sum operators
$\tP_k$ if $\max(d(1/p-1),1) <s<1/p$. This implies that the Haar system on $\mathbb{R}^d$ is not a basic sequence 
in $B_{p,q}^s(\mathbb{R}^d)$ for this parameter range, and strengthens the result of Triebel \cite{Tr1978}. 
As far as we know, the uniform boundedness
of the partial sum operators $\tP_k$ on $B_{p,q}^1(\mathbb{R}^d)$ has not been settled. This question is also open for $B_{p,q}^1(I^d)$, as
we only showed that $f-P_kf$ does not converge to zero
in the $B_{p,q}^1(I^d)$ quasi-norm for some $f\in B_{p,q}^1(I^d)$ but did not provide upper bounds for $s=1$).

\subsection{The limiting case $s=d(1/p-1)$}\label{sec33}
\subsubsection{Counterexamples for $p<q<\infty$}
Throughout this subsection, we fix $s=d(1/p-1)$, $d/(d+1)\le p<1$, and
$p<q<\infty$. In particular, this implies $s\le 1$ (with equality for
$p=d/(d+1)$).

We start with the statement in case b) in part ii) of Theorem \ref{theo2}.
The corresponding counterexamples have been suggested to me by T. Ullrich. They are similar to the counterexamples for Theorem \ref{theo1} but are now defined by linear combinations of special atoms. For the latter, we fix a function $a(x)=\prod_{i=1}^d \phi(x_i)$, where $\phi(x)\in C^\infty(\mathbb{R})$ is a univariate odd function, supported in $[-1,1]$, positive for $x\in (0,1)$, and such that $\|a\|_{\mathcal{C}^\sigma}=1$ for some $\sigma>1$. Obviously, if we define the functions $a_{j,\mathbf{i}}$, $\mathbf{i}\in \mathbb{Z}^d$ for $j\ge 1$ as in (\ref{Atoms}) from this $a$, then, with $C_0$ suitably fixed, they represent atoms of level $j\ge 1$, and we can apply Lemma \ref{lem2} for any $s\le 1$ to estimate $B_{p,q}^s(\mathbb{R}^d)$ quasi-norms of their linear combinations.

Consider the family of functions 
$$
g_k(x)=\sum_{j=1}^{n_k} j^{-1/p} 2^{(k+j)d} a_{k+j,\mathbf{i}_{k,j}}(x), \qquad n_k=(2^{k-2}-1)^d,\quad k=3,4,\ldots,
$$
where the multi-indices $\mathbf{i}_{k,j}$, $j=1,\ldots,n_k$, are chosen such that the support centers $x_{k,j}:=x_{k+j,\mathbf{i}_{k,j}}$ of the atoms $a_{k+j,\mathbf{i}_{k,j}}(x)$ are different, and coincide with the $n_k$ interior vertices of $T^d_{k-2}$. Note that $n_k\approx 2^{kd}$ as $k\to \infty$, and that $g_k(x)$ is a finite linear combination of atoms with different scale parameters whose supports are well-separated. 
By (\ref{BId}) and Lemma  \ref{lem2} we have $g_k(x)\in B_{p,q}^{d(1/p-1)}(I^d)$, $p<q<\infty$, with uniformly bounded quasi-norm for all $k\ge 3$ since
\bea
\|g_k\|_{B_{p,q}^{d(1/p-1)}}^q &\le& \|g_k\|_{B_{p,q}^{d(1/p-1)}(\mathbb{R}^d)}^q \le C(\|g_k\|^+_{B_{p,q}^{d(1/p-1)}(\mathbb{R}^d)})^q\\ 
&\le& C \sum_{j=1}^{n_k} 2^{(k+j)(d(1/p-1)-d/p)}(j^{-1/p} 2^{(k+j)d})^q\\
&=&C\sum_{j=1}^{n_k} j^{-q/p} \le C\sum_{j=1}^{n_k} j^{-q/p}<\infty.
\eea
For convenience, we use the same notation $g_k$ for the extension by zero of $g_k$ to $\mathbb{R}^d$.

On the other hand, since the centers $x_{k,j}$ of the atoms $a_{k+j,\mathbf{i}_{k,j}}$ are located at the interior vertices of $T_{k-2}^d$ and have supports in cubes of side-length $2^{-(k+j-1)}$, their supports are well-separated. Moreover, they have the same symmetry properties with respect to their centers as the function
$$
h_0^d(x):=h_0(x_1)\cdot \ldots \cdot h_0(x_d)
$$ 
has with respect to the origin. Here, $h_0$ is the univariate centralized Haar wavelet defined in (\ref{HW0}). Therefore, the Haar projection $P_kg_k$ onto
$S_k^d$ can easily be computed in explicit form: 
$$
(P_kg_k)(x)=\sum_{j=1}^{n_k} b_{k,j} j^{-1/p} 2^{(k+j)d} h_0^d(2^{k}(x-x_{k,j}),\quad k=3,4,\ldots,
$$
where $b_{k,j}$ is the average value of $a_{k+j,\mathbf{i}_{k,j}}$ over the cube in $T_k^d$ whose lowest vertex coincides with $x_{k,j}$.
This average value can easily be computed as 
$$
b_{k,j}=2^{kd} \int_{[0,2^{-k}]^d} a(2^{k+j}y)\,dy = 2^{kd} 2^{-(k+j)d}b_0 =b_02^{-jd},\qquad j=1,\ldots,n_k,
$$
where $b_0=(\int_0^1 \phi(x)\,dx)^d>0$ is a fixed constant.
Thus, the formula for $P_kg_k$ simplifies to
\be\label{Pkgk}
(P_kg_k)(x)=b_0 2^{kd}\sum_{j=1}^{n_k} j^{-1/p} h_0^d(2^{k}(x-x_{k,j})),\quad k=3,4,\ldots.
\ee

In order to get a lower bound for the  $B_{p,q}^{d(1/p-1)}$ quasi-norm of $P_kg_k$, we next compute a lower bound for its local mean
$\kappa(2^{-k},P_kg_k)(x)$, where the kernel $\kappa$ has the properties required for Lemma \ref{lem3} to hold. E.g., we could set $\kappa=a(2\cdot)$ with the above function $a$ because $s<1$ and the dilation factor $2$ ensures 
that $\mathrm{supp}(\kappa)\subset [-1/2,1/2]^d$. 
Since the support cubes (denoted by $I_{k,j}$) of the terms $h_0^d(2^{k}(x-x_{k,j}))$ in the representation (\ref{Pkgk}) have side-length $2^{-k+1}$ and are centered at $x_{k,j}\in T^d_{k-2}$, they are still well-separated, and we have
\bea
\kappa(2^{-k},P_kg_k)(x) &=& 2^{(k-1)d} j^{-1/p} (\kappa^{2^{-k}} \ast h_0^d(2^{k}(\cdot-x_{k,j})))(x)\\
&=& 2^{(k-1)d} j^{-1/p} 2^{kd}(\kappa(2^j\cdot) \ast h_0^d(2^{k}(\cdot))(x-x_{k,j})\\
&=& 2^{(k-1)d} j^{-1/p} (\kappa\ast h_0^d)(2^k(x-x_{k,j}),\qquad x\in I_{k,j},
\eea
where we have used (\ref{Atoms}). Since the cubes $I_{k,j}$ are also well-separated from the boundary of $I^d$, this lower bound holds for any extension $\tG_k\in S'(\mathbb{R}^d)$ of $P_kg_k$. Thus, since the $C^\infty$ function $\kappa\ast h_0^d$ is non-vanishing in a neighborhood of the origin by construction, we obtain
$$
\|\kappa(2^{-k},\tG_k)\|_{L_p(\mathbb{R}^d)}^p \ge c2^{(k-1)dp}\sum_{j=1}^{n_k} j^{-1} 2^{-kd} \ge c 2^{kd(p-1)}\log (n_k) \ge c2^{kd(p-1)} k.
$$  
By definition of Besov quasi-norms on domains and by Lemma \ref{lem2} we arrive at
$$
\|P_kg_k\|_{B_{p,q}^{d(1/p-1)}}
\ge c\inf_{\tG_k}  \|\tG_k\|^\ast_{B_{p,q}^{d(1/p-1)}(\mathbb{R}^d)} \\ 
\ge c 2^{kd(1/p-1)} \|\kappa(2^{-k},G_k)\|_{L_p(\mathbb{R}^d)}\ge c k^{1/p}.
$$
This shows that for $p<q\le 1$ the partial sum operators $P_k$, which are well-defined on the set of dyadic step-functions $S^d$ and  extend by continuity to $B_{p,q}^{d(1/p-1)}(I^d)$ due to the continuous embedding $B_{p,q}^{d(1/p-1)}(I^d)\subset L_1(I^d)$, cannot be uniformly bounded on $B_{p,q}^{d(1/p-1)}(I^d)$. This contradicts the Schauder basis property for $s=d(1/p-1)$, $p<q\le 1$,
and finishes the argument for case b).

For case a)
of part ii) of Theorem \ref{theo2}, we provide examples analogous to (\ref{gka}) in subsection \ref{sec23} which show that the average functionals $av_\Delta$ defined by (\ref{Av}) on the set of dyadic step functions $S^d$ cannot be extended to bounded linear functionals on
$B_{p,q}^{d(1/p-1)}(I^d)$ if $q>1$. For simplicity, consider $\Delta=I^d$,
and define $g_k=\tg_k|_{I^d}$ by the atomic decomposition
$$
\tg_k=\sum_{j=1}^k b_j a_{j,\mathbf{0}},   \qquad b_j=2^{jd}j^{-1}
$$
This is, up to different coefficient notation and the replacement of
characteristic functions $\chi_{\Delta_{j,\mathbf{2}}}$ by atoms
$a_{j,\mathbf{0}}$ defined in (\ref{Atoms}) with the above function $a$,
the same construction as in (\ref{gka}). Obviously, by construction
$$
av_{I^d}(g_k)=\sum_{j=1}^k b_j av_{I^d}(a(2^j\cdot))=
\sum_{j=1}^k b_j 2^{-jd}b_0 =b_0\sum_{j=1}^k j^{-1}\ge c\log(k)\to \infty
$$
as $k\to \infty$, while 
\bea
\|g_k\|_{B_{p,q}^{d(1/p-1)}}^q &\le& \|g_k\|_{B_{p,q}^{d(1/p-1)}(\mathbb{R}^d)}^q \le C(\|g_k\|^+_{B_{p,q}^{d(1/p-1)}(\mathbb{R}^d)})^q\\ 
&\le& C \sum_{j=1}^{k} (2^{j(d(1/p-1)-d/p)}j^{-1} 2^{jd})^q=C\sum_{j=1}^{n_k} j^{-q} \le C <\infty,
\eea
since $1<q<\infty$.

\subsubsection{Proof of the Schauder basis property for $0<q\le p$}   
We now turn to case a) in part ii) of Theorem \ref{theo2}, where 
$d/(d+1)<p<1$, $0<q\le p$, and $0<s=d(1/p-1) < 1$ can be assumed.
Since for this parameter range the set of all dyadic step functions is dense in $B_{p,q}^s(I^d)$,  it suffices to prove the uniform boundedness of the partial sum operators 
$$
\tilde{P}=\tilde{P}_k + \sum_{h\in \bar{H}_{k+1}^d} c_h(\cdot)h
$$ 
using Lemma \ref{lem4}, i.e., to establish (\ref{BTP}) for all $\tf \in Y\cap B_{p,q}^{d(1/p-1)}(\mathbb{R}^d)$. 
We proceed in several steps. 

\smallskip
\emph{Step 1}. Using the properties of $\tP\tf$, and in particular (\ref{AvTP}) and (\ref{AvH}), we show that
\be\label{Step1}
\sum_{j=k}^{\infty} 2^{jd(1/p-1)q}\|\kappa(2^{-j},\tP\tf)\|_{L_p(\mathbb{R}^d)}^q\le C s_{k+1}^q,\qquad 
\ee
where $s_{k+1}$ is given by
\be\label{Sk1}
s_{k+1}:=2^{-kd}\left(\sum_{\Delta',\Delta''\in \tT_{k+1}^d:
\mathrm{dist}_\infty(\Delta',\Delta'')=0} |av_{\Delta'}(\tf)-av_{\Delta''}(\tf)|^p\right)^{1/p}.
\ee

In the case $j\ge k$, consider any cube $\Delta\in \tT_k^d$, and denote the set of its neighbors in $T_k^d$ by
$$
n(\Delta)=\{\tilde{\Delta}\in T_k^d:\;\Delta\cap\tilde{\Delta}\neq\emptyset\}.
$$
Recall from (\ref{AvTP}) that
$\tP_k\tf|_\Delta=av_\Delta(\tf)$. Since $\kappa^{2^{-j}}(x-\cdot)$ is supported in a cube of side-length $2^{-j}$ centered at $x$ and is orthogonal to constants due to the assumed moment condition for the kernel $\kappa$, for $j\ge k$ and $x\in \Delta$ we have
\bea
|\kappa(2^{-j},\tP\tf)(x)|&=&|\kappa(2^{-j},\tP\tf(\cdot)-av_\Delta(\tf))(x)|\\
&\le& |\kappa(2^{-j},\tP_k\tf(\cdot)-av_\Delta(\tf))(x)|+
|\kappa(2^{-j},\tP\tf-\tP_k\tf)(x)|.
\eea
Here, both terms in the right-hand side vanish only if the support cube
of $\kappa^{2^{-j}}(x-\cdot)$ intersects with the boundary of any of the
dyadic cubes in $\tT_{k+1}^d$, where the piecewise constant functions $\tP\tf$ , $\tP_k\tf$ may have jumps. The set of these $x\in \Delta$ has measure
$\le C2^{-j}2^{-k(d-1)}$, and due to (\ref{AvTP}) 
 we have the bound
$$
|\kappa(2^{-j},\tP_k\tf(\cdot)-av_\Delta(\tf))(x)|\le C\sum_{\tilde{\Delta}\in n(\Delta)} |av_\Delta(\tf)-av_{\tilde{\Delta}}(\tf)|,
$$
which in turn can be estimated by the sum of differences $
|av_{\Delta'}(\tf)-va_{\Delta''}(\tf)|$ appearing in (\ref{Step1}) with neighboring $\Delta',\Delta''\in \tT_{k+1}^d$ belonging to the union of cubes in
$n(\Delta)$. The other term is similarly bounded 
since for $x\in \Delta$
$$
|\kappa(2^{-j},\tP\tf-\tP_k\tf)(x)|\le C
\| \sum_{h\in \bar{H}_{k+1}^d: } |c_h(\tf)h(\cdot)| \|_{L_\infty(\cup_{\tilde{\Delta}\in n(\Delta)}\tilde{\Delta})},
$$
and we can apply (\ref{AvH}).
Thus, altogether we arrive at
\be\label{Step1a}
|\kappa(2^{-j},\tP\tf)(x)|\le C {\sum}'_{\Delta',\Delta''} |av_{\Delta'}(\tf)-av_{\Delta''}(\tf)|,
\qquad x\in \Delta\in \tT_k^d, \qquad j\ge k,
\ee
where $\sum'$ indicates that the summation extends to all those neighboring
dyadic cubes $\Delta',\Delta''$ in $\tT_{k+1}^d$ which belong to  the union of all cubes in
$n(\Delta)$. This bound is only needed on a subset of $\Delta$ of
measure $\le C2^{-j}2^{-k(d-1)}$.
 
From (\ref{Step1a}) we get for $j\ge k$
\bea
\|\kappa(2^{-j},\tP\tf)\|_{L_p(\mathbb{R}^d)}^p&=& \sum_{\Delta\in \tT^d_{k}}
\|\kappa(2^{-j},\tP\tf)\|_{L_p(\Delta)}^p\\
&\le&C \sum_{\Delta\in \tT^d_{k}} 2^{-j-k(d-1)} 
{\sum}'_{\Delta',\Delta''} |av_{\Delta'}(\tf)-av_{\Delta''}(\tf)|^p\\
&\le&C 2^{k-j}2^{-kd(1-p)} \sum_{\Delta',\Delta''\in \tT^d_{k+1}:\,\mathrm{dist}_\infty(\Delta',\Delta'')=0} |av_{\Delta'}(\tf)-av_{\Delta''}(\tf)|^p\\
&=& C2^{-j}2^{-k(d-1)}s_{k+1}^p,
\eea
where we have used that each term $|av_{\Delta'}(\tf)-av_{\Delta''}(\tf)|^p$  belongs to at most $3^d$ neighborhoods $n(\Delta)$.
Taking the previous estimate to the power $q/p$ and substituting the result into the left-hand side of (\ref{Step1}) leads to the desired estimate in (\ref{Step1}). Indeed,
$$
\sum_{j=k}^{\infty} 2^{jd(1/p-1)q}\|\kappa(2^{-j},\tP\tf)\|_{L_p(\mathbb{R}^d)}^q\le C2^{k(1/p-d/p)q}s_{k+1}^q \sum_{j=k}^\infty 2^{-j(1/p-d(1/p-1))q}
\le C2^{-kdq}s_{k+1}^q,
$$
since $1/p-d(1/p-1)<0$ for $d/(d+1)<p<1$.

\smallskip
\emph{Step 2}. For $1\le j<k$, we start with
$$
\|\kappa(2^{-j},\tP\tf)\|_{L_p(\mathbb{R}^d)}^p\le
\|\kappa(2^{-j},\tf-\tP\tf)\|_{L_p(\mathbb{R}^d)}^p+\|\kappa(2^{-j},\tf)\|_{L_p(\mathbb{R}^d)}^p,
$$
(similarly for $j=0$ and $\kappa_0(1,\tP\tf)$), and proceed with estimates for the term corresponding to $\tf-\tP\tf$ (after substitution into the expression for the local means quasi-norm, the other term will be automatically bounded by
the right-hand side in (\ref{BTP})).  This time we use the fact that 
$\tf-\tP\tf$ has zero average on each dyadic cube $\tilde{\Delta}\in \tT_k^d$,
and that the kernel $\kappa$ is smooth. With the short-hand notation 
$\Delta_x^j$ for the support cube of $\kappa^{2^{-j}}(x-\cdot)$, this yields
\bea
|\kappa(2^{-j},\tf-\tP\tf)(x)|&\le& \sum_{\tilde{\Delta}\in \tT^d_k: 
\,\tilde{\Delta}\cap \Delta_x^j\neq \emptyset} \left|\int_{\tilde{\Delta}} 2^{jd} \kappa^{2^{-j}}(x-y)(\tf-\tP\tf)(y)\,dy\right| \\
&\le& \sum_{\tilde{\Delta}\in \tT^d_{k}:\tilde{\Delta}\cap \Delta_x^j\neq \emptyset} 
\inf_{\xi\in\mathbb{R}}\|\kappa^{2^{-j}}(x-\cdot)-\xi\|_{L_\infty(\Delta_x^j\cap \tilde{\Delta})}\int_{\tilde{\Delta}} |(\tf-\tP\tf)(y)|\,dy\\
&\le& C 2^{jd}2^{j-k} \sum_{\tilde{\Delta}\in \tT^d_{k}:
\tilde{\Delta}\cap \Delta_x^j\neq \emptyset} 
\int_{\tilde{\Delta}} |\tf-\tP\tf)(y)|\,dy\\
&\le& C 2^{jd}2^{j-k} \sum_{\tilde{\Delta}\in \tT^d_{k}:\, \mathrm{dist}_\infty(\tilde{\Delta},\Delta)\le C_1 2^{-j}} \int_{\tilde{\Delta}} |\tf-\tP\tf|\,dy,\qquad x\in \Delta\in \tT_k^d,
\eea
if the constant $C_1$ is suitably chosen depending on $d$. 
Since moment conditions of the kernel $\kappa$ did not play a role in this part, the estimate will also hold for $j=0$ and $\kappa_0(1,\tf-\tP\tf)$.

We next compute the $L_p(\mathbb{R}^d)$ quasi-norm of $\kappa(2^{-j},\tf-\tP\tf)$:
\bea
&&\|\kappa(2^{-j},\tf-\tP\tf)\|_{L_p(\mathbb{R}^d)}^p =
\sum_{\Delta\in \tT_k^d} \|\kappa(2^{-j},\tf-\tP\tf)\|_{L_p(\Delta)}^p\\
&& \qquad\le C2^{jdp}2^{(j-k)p} \sum_{\Delta\in \tT_k^d} 2^{-kd} \left(\sum_{\tilde{\Delta}\in \tT^d_{k}:\, \mathrm{dist}_\infty(\tilde{\Delta},\Delta)\le C_1 2^{-j}} \int_{\tilde{\Delta}} |\tf-\tP\tf|\,dy\right)^p\\
&& \qquad\le C2^{jdp}2^{(j-k)p} \sum_{\Delta\in \tT_k^d} 2^{-kd} \sum_{\tilde{\Delta}\in \tT^d_{k}:\, \mathrm{dist}_\infty(\tilde{\Delta},\Delta)\le C_1 2^{-j}} (\int_{\tilde{\Delta}} |\tf-\tP\tf|\,dy)^p\\
&& \qquad\le C2^{jd(p-1)}2^{(j-k)p}\sum_{\tilde{\Delta}\in \tT^d_{k}} (\int_{\tilde{\Delta}} |\tf-\tP\tf|\,dy)^p\\ 
&=& \qquad C2^{jd(p-1)}2^{(j-k)p} \bar{s}_{k}^p,  \qquad j=1,\ldots,k-1,
\eea
where in the change of summation step we used that the number appearances of integrals over any fixed $\tilde{\Delta}\in \tT_{k}^d$
is bounded by $C2^{d(k-j)}$. The notation
\be\label{BarSk}
\bar{s}_{k}:=
\left(\sum_{\tilde{\Delta}\in \tT^d_{k}} (\int_{\tilde{\Delta}} |\tf-\tP\tf|\,dy)^p\right)^{1/p}
\ee
is introduced for convenience. The estimate also holds for $j=0$ with $\kappa$ replaced by $\kappa_0$.

This eventually gives
$$
\|\kappa_0(1,\tP\tf)\|_{L_p(\mathbb{R}^d)}^q +\sum_{j=1}^{k-1} 2^{jd(1/p-1)q}\|\kappa(2^{-j},\tP\tf)\|_{L_p(\mathbb{R}^d)}^q   \le C2^{-kq} \bar{s}_{k}^q \sum_{j=0}^{k-1} 2^{jq} \le C \bar{s}_{k}^q.
$$
Together with (\ref{Step1}), we arrive at
\be\label{Step2}
\|\tP\tf\|^\ast_{B^{d(1/p-1)}_{p,q}(\mathbb{R}^d)} \le C(s_{k+1} +\bar{s}_{k}+\|\tf\|^\ast_{B^{d(1/p-1)}_{p,q}(\mathbb{R}^d)}).
\ee

\smallskip
\emph{Step 3}. It remains to deal with the terms $s_{k+1}$ and $\bar{s}_{k}$ in (\ref{Step2}) which do not depend on $q$. This task is reminiscent of 
the estimation of the right-hand side in (\ref{Step11}) in the proof of Theorem \ref{theo1}. We show all details for $\bar{s}_{k}$, the estimates for $s_{k+1}$ are analogous, we only indicate the changes in the argument.

We explore the atomic decomposition (\ref{AtRep}) of $\tf\in Y\cap B^{d(1/p-1)}_{p,q}(\mathbb{R}^d)$, and observe that for $0<q\le p$ we have 
$$
\|\tf\|^+_{B^{d(1/p-1)}_{p,p}(\mathbb{R}^d)}\le C \|\tf\|^+_{B^{d(1/p-1)}_{p,q}(\mathbb{R}^d)}.
$$
Therefore, it suffices to set $q=p$ and to show that
\be\label{Step3}
\bar{s}_k^p \le C \sum_{j=0}^\infty 2^{-jdp} c_j^p,\qquad
c_j:=\left(\sum_{\mathbf{i}\in \mathbb{Z}^d }
|c_{j,\mathbf{i}}|^p\right)^{1/p},
\ee
since, after taking the infimum in (\ref{Step3}) with respect to all atomic decompositions representing the same $\tf$, we get the desired bound 
$$
\bar{s}_k\le C\|\tf\|^+_{B_{p,p}^{d(1/p-1)}(\mathbb{R}^d)}\le C \|\tf\|^+_{B^{d(1/p-1)}_{p,q}(\mathbb{R}^d)},\qquad 0<q\le p.
$$

For each integral over a dyadic cube $\tilde{\Delta}\in\tT_{k}^d$ in $\bar{s}_k^p$, we estimate
$$
\int_{\tilde{\Delta}} |\tf-\tP\tf|\,dx  \le \sum_{j=0}^\infty 
\sum_{\mathbf{i}:\,\mathrm{supp}(a_{j,\mathbf{i}})\cap \tilde{\Delta}\neq\emptyset} |c_{j,\mathbf{i}}| \int_{\tilde{\Delta}} |a_{j,\mathbf{i}}-\tP a_{j,\mathbf{i}}|\,dx.
$$
For $j> k$, we can estimate the relevant terms in the sum by
$$
\int_{\tilde{\Delta}} |a_{j,\mathbf{i}}-\tP a_{j,\mathbf{i}}| \,dx\le C 
\int_{\tilde{\Delta}} |a_{j,\mathbf{i}}| \,dx \le C 2^{-jd},
$$
and each such term may appear only for $\le C$ different $\tilde{\Delta}$ (this $C$ depends on $C_0$ and $d$). For $j\le k$, we explore  the $C^1$
continuity of the atoms (recall that we assumed $a\in \mathcal{C}^\sigma_{C_0}$ with $\sigma>1$ in (\ref{Atoms})) which gives
\bea
|a_{j,\mathbf{i}}(x)-\tP a_{j,\mathbf{i}}(x)|&\le& 
|a_{j,\mathbf{i}}(x)- av_{\tilde{\Delta}}(a_{j,\mathbf{i}})| +
\sum_{h\in \bar{H}_{k+1}^d:\,\mathrm{supp}(h)=\tilde{\Delta}} 2^{kd/2} |c_h(a_{j,\mathbf{i}})|\\
&\le& C 2^{j-k},\qquad x\in \tilde{\Delta},
\eea
where $C$ depends on $C_0, \sigma$, and $d$. Thus, in this case we get
$$
\int_{\tilde{\Delta}} |a_{j,\mathbf{i}}-\tP a_{j,\mathbf{i}}| \,dx\le C 2^{j-k}2^{-kd}, 
$$
where each such term appears for $\le C2^{(k-j)d}$ different $\tilde{\Delta}$.

Substitution into the expression (\ref{BarSk}) for $\bar{s}_k^p$ in  gives
\bea
\bar{s}_k^p&=&\sum_{\tilde{\Delta}\in \tT^d_{k}} (\int_{\tilde{\Delta}} |\tf-\tP\tf|\,dy)^p\\
&\le& C\left(\sum_{j=0}^k \sum_{\mathbf{i}\in \mathbb{Z}^d} 2^{(k-j)d} 2^{(j-k-kd)p}
|c_{j,\mathbf{i}}|^p + \sum_{j=k+1}^\infty \sum_{\mathbf{i}\in \mathbb{Z}^d} 2^{-jdp} |c_{j,\mathbf{i}}|^p\right)\\
&=& C\left(2^{k(-1+d(1/p-1))p} \sum_{j=0}^k 2^{j(-d/p+1)p} c_j^p + \sum_{j=k+1}^\infty 2^{-jdp}c_j^p\right).
\eea
Since $1-d(1/p-1)<0$ for our parameter range $p>d/(d+1)$ implies 
$$
2^{k(-1+d(1/p-1))p}2^{j(-d/p+1)p}=2^{-(k-j)(1-d(1/p-1))p}2^{-jdp}\le 2^{-jdp},
$$
this proves (\ref{Step3}).

To estimate $s_{k+1}^p$ by the right-hand side in (\ref{Step3}), instead of the terms $\|\tf-\tP\tf\|_{L_1(\tilde{\Delta})}^p$ with $\tilde{\Delta}\in \tT_k^d$, we must now consider the terms
\bea
&&2^{-kdp}|av_{\Delta'}(\tf)-av_{\Delta''}(\tf)|^p=2^{dp}
\left|\int_{\Delta'}  \tf\,dy -\int_{\Delta''}\tf \,dx\right|^p\\
&&\qquad \le 2^{dp}\sum_{j=0}^\infty \sum_{\mathbf{i}\in\mathbb{Z}^d:\,
\mathrm{supp}(a_{j,\mathbf{i}}) \cap(\Delta'\cup\Delta'') \neq\emptyset}
|c_{j,\mathbf{i}}|^p\left|\int_{\Delta'}  a_{j,\mathbf{i}}\,dy -\int_{\Delta''} a_{j,\mathbf{i}} \,dx\right|^p
\eea
for neighboring dyadic cubes $\Delta',\Delta''$ in $\tT_{k+1}$. 
But the estimates for the quantities
$$
\left|\int_{\Delta'}  a_{j,\mathbf{i}}\,dy -\int_{\Delta''} a_{j,\mathbf{i}}\,dx\right|^p
$$
in the two cases $j>k$ and $j\le k$ look completely the same as the estimates
for $\|a_{j,\mathbf{i}}-\tP a_{j,\mathbf{i}}\|_{L_1(\tilde{\Delta})}$.
The remaining steps can be repeated without change. 

Together with (\ref{Step2}), we have shown the uniform boundedness of
the partial sum operators $\tP$ in $B_{p,q}^s(\mathbb{R}^d)$. This finishes the argument.

\end{document}